\documentclass[reqno,12pt,a4paper]{amsart}

\usepackage{mathrsfs}
\usepackage{amssymb}
\usepackage{amsfonts}
\usepackage{latexsym}
\usepackage{amsthm}

\usepackage{graphicx}
\def\lb{\label}

\newcommand{\er}[1]{\textrm{(\ref{#1})}}

\begin{document}

%%%%%%%%%% Some definitions %%%%%%%%%%

%%%%%%%% Equations, theorems %%%%%%%%%
\renewcommand{\theequation}{\arabic{section}.\arabic{equation}}
\theoremstyle{plain}
\newtheorem{theorem}{\bf Theorem}[section]
\newtheorem{lemma}[theorem]{\bf Lemma}
\newtheorem{corollary}[theorem]{\bf Corollary}
\newtheorem{proposition}[theorem]{\bf Proposition}
\newtheorem{definition}[theorem]{\bf Definition}
\newtheorem{remark}[theorem]{\it Remark}
%\theoremstyle{remark}
%\newtheorem{remark}[theorem]{\bf Remark}

%%%%% Alphabet %%%%%
\def\a{\alpha}  \def\cA{{\mathcal A}}     \def\bA{{\bf A}}  \def\mA{{\mathscr A}}
\def\b{\beta}   \def\cB{{\mathcal B}}     \def\bB{{\bf B}}  \def\mB{{\mathscr B}}
\def\g{\gamma}  \def\cC{{\mathcal C}}     \def\bC{{\bf C}}  \def\mC{{\mathscr C}}
\def\G{\Gamma}  \def\cD{{\mathcal D}}     \def\bD{{\bf D}}  \def\mD{{\mathscr D}}
\def\d{\delta}  \def\cE{{\mathcal E}}     \def\bE{{\bf E}}  \def\mE{{\mathscr E}}
\def\D{\Delta}  \def\cF{{\mathcal F}}     \def\bF{{\bf F}}  \def\mF{{\mathscr F}}
\def\c{\chi}    \def\cG{{\mathcal G}}     \def\bG{{\bf G}}  \def\mG{{\mathscr G}}
\def\z{\zeta}   \def\cH{{\mathcal H}}     \def\bH{{\bf H}}  \def\mH{{\mathscr H}}
\def\e{\eta}    \def\cI{{\mathcal I}}     \def\bI{{\bf I}}  \def\mI{{\mathscr I}}
\def\p{\psi}    \def\cJ{{\mathcal J}}     \def\bJ{{\bf J}}  \def\mJ{{\mathscr J}}
\def\vT{\Theta} \def\cK{{\mathcal K}}     \def\bK{{\bf K}}  \def\mK{{\mathscr K}}
\def\k{\kappa}  \def\cL{{\mathcal L}}     \def\bL{{\bf L}}  \def\mL{{\mathscr L}}
\def\l{\lambda} \def\cM{{\mathcal M}}     \def\bM{{\bf M}}  \def\mM{{\mathscr M}}
\def\L{\Lambda} \def\cN{{\mathcal N}}     \def\bN{{\bf N}}  \def\mN{{\mathscr N}}
\def\m{\mu}     \def\cO{{\mathcal O}}     \def\bO{{\bf O}}  \def\mO{{\mathscr O}}
\def\n{\nu}     \def\cP{{\mathcal P}}     \def\bP{{\bf P}}  \def\mP{{\mathscr P}}
\def\r{\rho}    \def\cQ{{\mathcal Q}}     \def\bQ{{\bf Q}}  \def\mQ{{\mathscr Q}}
\def\s{\sigma}  \def\cR{{\mathcal R}}     \def\bR{{\bf R}}  \def\mR{{\mathscr R}}
\def\S{\Sigma}  \def\cS{{\mathcal S}}     \def\bS{{\bf S}}  \def\mS{{\mathscr S}}
\def\t{\tau}    \def\cT{{\mathcal T}}     \def\bT{{\bf T}}  \def\mT{{\mathscr T}}
\def\f{\phi}    \def\cU{{\mathcal U}}     \def\bU{{\bf U}}  \def\mU{{\mathscr U}}
\def\F{\Phi}    \def\cV{{\mathcal V}}     \def\bV{{\bf V}}  \def\mV{{\mathscr V}}
\def\P{\Psi}    \def\cW{{\mathcal W}}     \def\bW{{\bf W}}  \def\mW{{\mathscr W}}
\def\o{\omega}  \def\cX{{\mathcal X}}     \def\bX{{\bf X}}  \def\mX{{\mathscr X}}
\def\x{\xi}     \def\cY{{\mathcal Y}}     \def\bY{{\bf Y}}  \def\mY{{\mathscr Y}}
\def\X{\Xi}     \def\cZ{{\mathcal Z}}     \def\bZ{{\bf Z}}  \def\mZ{{\mathscr Z}}
\def\O{\Omega}

\newcommand{\gA}{\mathfrak{A}}
\newcommand{\gB}{\mathfrak{B}}
\newcommand{\gC}{\mathfrak{C}}
\newcommand{\gD}{\mathfrak{D}}
\newcommand{\gE}{\mathfrak{E}}
\newcommand{\gF}{\mathfrak{F}}
\newcommand{\gG}{\mathfrak{G}}
\newcommand{\gH}{\mathfrak{H}}
\newcommand{\gI}{\mathfrak{I}}
\newcommand{\gJ}{\mathfrak{J}}
\newcommand{\gK}{\mathfrak{K}}
\newcommand{\gL}{\mathfrak{L}}
\newcommand{\gM}{\mathfrak{M}}
\newcommand{\gN}{\mathfrak{N}}
\newcommand{\gO}{\mathfrak{O}}
\newcommand{\gP}{\mathfrak{P}}
\newcommand{\gR}{\mathfrak{R}}
\newcommand{\gS}{\mathfrak{S}}
\newcommand{\gT}{\mathfrak{T}}
\newcommand{\gU}{\mathfrak{U}}
\newcommand{\gV}{\mathfrak{V}}
\newcommand{\gW}{\mathfrak{W}}
\newcommand{\gX}{\mathfrak{X}}
\newcommand{\gY}{\mathfrak{Y}}
\newcommand{\gZ}{\mathfrak{Z}}

\def\ve{\varepsilon}   \def\vt{\vartheta}    \def\vp{\varphi}    \def\vk{\varkappa}

\def\Z{{\mathbb Z}}    \def\R{{\mathbb R}}   \def\C{{\mathbb C}}
\def\T{{\mathbb T}}    \def\N{{\mathbb N}}   \def\dD{{\mathbb D}}

%%%%% Arrows %%%%%

\def\la{\leftarrow}              \def\ra{\rightarrow}            \def\Ra{\Rightarrow}
\def\ua{\uparrow}                \def\da{\downarrow}
\def\lra{\leftrightarrow}        \def\Lra{\Leftrightarrow}

%%%%% Typography %%%%%

\def\lt{\biggl}                  \def\rt{\biggr}
\def\ol{\overline}               \def\wt{\widetilde}
\def\no{\noindent}

%%%%% Math signs %%%%%

\let\ge\geqslant                 \let\le\leqslant
\def\lan{\langle}                \def\ran{\rangle}
\def\/{\over}                    \def\iy{\infty}
\def\sm{\setminus}               \def\es{\emptyset}
\def\ss{\subset}                 \def\ts{\times}
\def\pa{\partial}                \def\os{\oplus}
\def\om{\ominus}                 \def\ev{\equiv}
\def\iint{\int\!\!\!\int}        \def\iintt{\mathop{\int\!\!\int\!\!\dots\!\!\int}\limits}
\def\el2{\ell^{\,2}}             \def\1{1\!\!1}
\def\sh{\sharp}
\def\wh{\widehat}
%%%%% Math operations %%%%%

\def\where{\mathop{\mathrm{where}}\nolimits}
\def\as{\mathop{\mathrm{as}}\nolimits}
\def\Area{\mathop{\mathrm{Area}}\nolimits}
\def\arg{\mathop{\mathrm{arg}}\nolimits}
\def\const{\mathop{\mathrm{const}}\nolimits}
\def\det{\mathop{\mathrm{det}}\nolimits}
\def\diag{\mathop{\mathrm{diag}}\nolimits}
\def\diam{\mathop{\mathrm{diam}}\nolimits}
\def\dim{\mathop{\mathrm{dim}}\nolimits}
\def\dist{\mathop{\mathrm{dist}}\nolimits}
\def\Im{\mathop{\mathrm{Im}}\nolimits}
\def\Iso{\mathop{\mathrm{Iso}}\nolimits}
\def\Ker{\mathop{\mathrm{Ker}}\nolimits}
\def\Lip{\mathop{\mathrm{Lip}}\nolimits}
\def\rank{\mathop{\mathrm{rank}}\limits}
\def\Ran{\mathop{\mathrm{Ran}}\nolimits}
\def\Re{\mathop{\mathrm{Re}}\nolimits}
\def\Res{\mathop{\mathrm{Res}}\nolimits}
\def\res{\mathop{\mathrm{res}}\limits}
\def\sign{\mathop{\mathrm{sign}}\nolimits}
\def\span{\mathop{\mathrm{span}}\nolimits}
\def\supp{\mathop{\mathrm{supp}}\nolimits}
\def\Tr{\mathop{\mathrm{Tr}}\nolimits}
\def\BBox{\hspace{1mm}\vrule height6pt width5.5pt depth0pt \hspace{6pt}}

%%%%%%%%%%%%% specialities %%%%%%%%%%%%%%

\newcommand\nh[2]{\widehat{#1}\vphantom{#1}^{(#2)}}
%{{\mathop{#1}\limits^\wedge}\vphantom{#1}^{(#2)}}
\def\dia{\diamond}

\def\Oplus{\bigoplus\nolimits}

%%%%%%%%%%% End of definitions %%%%%%%%%%

%%%%% OLD OLD OLD

\def\qqq{\qquad}
\def\qq{\quad}
\let\ge\geqslant
\let\le\leqslant
\let\geq\geqslant
\let\leq\leqslant
\newcommand{\ca}{\begin{cases}}
\newcommand{\ac}{\end{cases}}
\newcommand{\ma}{\begin{pmatrix}}
\newcommand{\am}{\end{pmatrix}}
\renewcommand{\[}{\begin{equation}}
\renewcommand{\]}{\end{equation}}
\def\bu{\bullet}

\title[{Schr\"odinger operator with periodic plus
compactly supported potentials}]
{1D  Schr\"odinger operator with periodic plus
compactly supported potentials}

\date{\today}
\author[Evgeny Korotyaev]{Evgeny Korotyaev}
\address{School of Math., Cardiff University.
Senghennydd Road, CF24 4AG Cardiff, Wales, UK.
email \ KorotyaevE@cf.ac.uk,
{\rm Partially supported by EPSRC grant EP/D054621.}}


\begin{abstract}
We consider the 1D Schr\"odinger operator $Hy=-y''+(p+q)y$ with a periodic potential $p$  plus
compactly supported potential $q$ on the real line.
The spectrum of $H$ consists of an absolutely continuous part plus a finite number of simple eigenvalues in each spectral gap $\g_n\ne \es, n\geq 0$, where $\g_0$ is unbounded gap. We prove the following results:
 1)  we determine the distribution of resonances in the disk with large radius,
 2) a forbidden domain for the resonances is specified,
3) the asymptotics of eigenvalues and antibound states are determined,
4) if $q_0=\int_\R qdx=0$, then roughly speaking in each nondegenerate gap $\g_n$ for $n$ large enough there are two eigenvalues
and zero antibound state or zero eigenvalues and two antibound states,
5) if $H$ has infinitely many gaps in the continuous spectrum, then for any sequence $\s=(\s)_1^\iy, \s_n\in \{0,2\}$, there exists a compactly supported potential $q$ such that $H$ has $\s_n$ bound states and $2-\s_n$ antibound states in each gap $\g_n$ for $n$ large enough.
6) For any $q$ (with $q_0=0$), $\s=(\s_n)_{1}^\iy$, where $\s_n\in \{0,2\}$ and   for any sequence  $\d=(\d_n)_1^\iy\in \ell^2, \d_n>0$
there exists a potential $p\in L^2(0,1)$  such that each gap length  $|\g_n|=\d_n, n\ge 1$ and $H$ has exactly $\s_n$ eigenvalues and $2-\s_n$ antibound states in each gap $\g_n\ne \es$ for $n$ large enough.

\end{abstract}

\maketitle

\section{Introduction}

Consider the Schr\"odinger operator $H=H_0+q$ acting in $L^2(\R)$, where $H_0=-{d^2\/dx^2}+p(x)$
and $p\in L^2(0,1)$ is the real 1-periodic potential. The  real compactly supported potential $q$ belongs to the
class $\cQ_{t}=\{q\in L^2(\R ): \  [0,t] \ is\ the\ convex\ hull\ of\ the\ support\ of\ q\}$ for some $t>0$. The spectrum of $H_0$ consists of spectral bands $\gS_n$ and is given by (see Fig. 1)
$$
\s(H_0)=\s_{ac}(H_0)=\cup \gS_n,\qq \gS_n=[E^+_{n-1},E^-_n],n\ge 1,\qq
E_0^+<..\le E^+_{n-1}< E^-_n \le E^+_{n}<...
$$
We assume that $E_0^+=0$.  The sequence $E_0^+<E_1^-\le E_1^+\ <\dots$  is the spectrum of the equation
\[
\lb{1}
-y''+p(x)y=\l y, \ \ \ \ \l\in \C ,
\]
with the 2-periodic boundary conditions, i.e. $y(x+2)=y(x), x\in \R$.
The bands $\gS_n, \gS_{n+1}$ are separated by  a gap $\g_{n}=(E^-_{n},E^+_n)$ and let $\g_0=(-\iy,E_0^+)$.
If a gap degenerates, that is $\g_n=\es $, then the corresponding bands $\gS_{n} $ and $\gS_{n+1}$ merge.
If $E_n^-=E_n^+$ for some $n$, then this number $E_n^{\pm}$ is the double eigenvalue of equation \er{1} with the 2-periodic boundary conditions. The lowest eigenvalue $E_0^+=0$ is always simple  and the corresponding eigenfunction is 1-periodic. The eigenfunctions, corresponding to the eigenvalue $E_{2n}^{\pm}$, are 1-periodic, and for the case $E_{2n+1}^{\pm}$ they are anti-periodic,  i.e., $y(x+1)=-y(x),\ \ x\in\R$.

It is well known, that the spectrum of $H$ consists of an absolutely continuous part $\s_{ac}(H)=\s_{ac}(H_0)$ plus a finite number of simple eigenvalues in each gap $\g_n\ne \es, n\geq 0$, see  \cite{Rb}, \cite{F1}.
Moreover, in every open gap $\g_n$ for $n$ large enough  the operator $H$ has  at most two bound states \cite{Rb} and precisely one bound state in the case $q_0=\int_\R q(x)dx\ne 0$ \cite{Zh}, \cite{F2}.
Note that  the potential $q$ in \cite{Rb},\cite{Zh} belongs to the more general class, see also \cite{GS}, \cite{So}.

\begin{figure}
\tiny
\unitlength=1.00mm
\special{em:linewidth 0.4pt}
\linethickness{0.4pt}
\begin{picture}(108.67,33.67)
%coordinate lines
\put(41.00,17.33){\line(1,0){67.67}}
\put(44.33,9.00){\line(0,1){24.67}}
\put(108.33,14.00){\makebox(0,0)[cc]{$\Re\l$}}
\put(41.66,33.67){\makebox(0,0)[cc]{$\Im\l$}}
\put(42.00,14.33){\makebox(0,0)[cc]{$0$}}
%spectrum
\put(44.33,17.33){\linethickness{4.0pt}\line(1,0){11.33}}
%\put(48.66,17.33){\linethickness{4.0pt}\line(1,0){11.33}}
\put(66.66,17.33){\linethickness{4.0pt}\line(1,0){11.67}}
\put(82.00,17.33){\linethickness{4.0pt}\line(1,0){12.00}}
\put(95.66,17.33){\linethickness{4.0pt}\line(1,0){11.00}}
%endpoints of gaps
\put(46.66,20.00){\makebox(0,0)[cc]{$E_0^+$}}
%\put(48.66,20.00){\makebox(0,0)[cc]{$E_0^+$}}
\put(56.66,20.33){\makebox(0,0)[cc]{$E_1^-$}}
%\put(59.66,20.33){\makebox(0,0)[cc]{$E_1^-$}}
\put(68.66,20.33){\makebox(0,0)[cc]{$E_1^+$}}
%\put(66.66,20.33){\makebox(0,0)[cc]{$E_1^+$}}
\put(78.33,20.33){\makebox(0,0)[cc]{$E_2^-$}}
\put(84.33,20.33){\makebox(0,0)[cc]{$E_2^+$}}
%\put(82.33,20.33){\makebox(0,0)[cc]{$E_2^+$}}
\put(93.00,20.33){\makebox(0,0)[cc]{$E_3^-$}}
\put(98.66,20.33){\makebox(0,0)[cc]{$E_3^+$}}
%\put(96.66,20.33){\makebox(0,0)[cc]{$E_3^+$}}
\put(106.33,20.33){\makebox(0,0)[cc]{$E_4^-$}}
\end{picture}
\caption{The cut domain $\C\sm \cup \gS_n$ and the cuts (bands) $\gS_n=[E^+_{n-1},E^-_n], n\ge 1$}
\lb{sS}
\end{figure}

Let $\vp(x,z), \vt(x,z)$ be the solutions of the equation $-y''+py=z^2y$
satisfying $\vp'(0,z)=\vt(0,z)=1$ and $\vp(0,z)=\vt'(0,z)=0$, where $u'=\pa_x u$.
The Lyapunov function is defined by $\D(z)={1\/2}(\vp'(1,z)+\vt(1,z))$.
The function $\D^2(\sqrt \l)$ is entire, where $\sqrt \l$ is defined by $\sqrt 1=1$. Introduce the function $\O(\l)=(1-\D^2(\sqrt \l))^{1\/2}, \l\in \ol\C_+$ and we fix the branch $\O(\l)=(1-\D^2(\sqrt \l))^{1\/2}$ by the condition $\O(\l+i0)>0$ for $\l\in \gS_1=[E^+_{0},E^-_1]$.
Introduce the two-sheeted Riemann surface $\L$ of $\O(\l)=(1-\D^2(\sqrt \l))^{1\/2}$ obtained by joining the upper and lower rims of two copies of the cut plane $\C\sm\s_{ac}(H_0)$ in the usual (crosswise) way.
The n-th gap on the first physical sheet $\L_1$ we will denote by $\g_n^{(1)}$ and the same gap but on the second
 nonphysical sheet $\L_2$ we will denote by $\g_n^{(2)}$ and let $\g_n^{(0)}$ be the union of $\ol\g_n^{(1)}$ and $\ol\g_n^{(2)}$, i.e.,
 $$
 \g_n^{(0)}=\ol\g_n^{(1)}\cup \ol\g_n^{(2)}.
 $$

Introduce the function $D(\l)=\det (I+q(H_0-\l)^{-1})$, which is meromorphic
in $\L$, see \cite{F1}. Recall that   $1/D$ is the transmission coefficient in the S-matrix for the operators $H,H_0$, see Sect. 2. If $D$ has some poles, then they coincide with some $E_n^\pm$. It is well known that if $D(\l)=0$ for some zero $\l\in \L_1$, then $\l$ is an eigenvalue of $H$ and $\l\in\cup \g_n^{(1)}$. Note that there are no eigenvalues on the spectrum $\s_{ac}(H_0)\ss \L_1$, since $|D(\l)|\ge 1$  on $\s_{ac}(H_0)\ss \L_1$ (see \er{iab} and all these facts in \cite{F1}).

Define the functions $A, J$ by
\[
J(\l)= 2\O(\l+i0)\Im D(\l+i0),\qqq \  A(\l)=\Re D(\l+i0)-1, \qq  for \qq \l\in \s(H_0)\ss\L_1.
\]
In Lemma \ref{T32} we will show that $A, J$ are entire functions on $\C$ and they are real on the real line. Instead of the function $D$ we will study the {\bf modified function} $\Xi=2i\O D$ on $\L$. We will show that $\Xi$ satisfies
\[
\lb{T1-1}
\Xi=2i\O D=2i \O (1+A)-J\qqq on \qqq  \L,
\]
see Theorem \ref{T1}.
Recall that $\O$ is analytic  in $\L$ and $\O(\l)=0$ for some $\l\in \L$ iff $\l=E_n^-$ or $\l=E_n^+$  for some $n\ge 0$.  Then the function $\Xi$ is analytic on $\L$ and has branch points $E_n^\pm, \g_n\ne \es$. The zeros of $\Xi$ define the eigenvalues and resonances, similar to the case $p=0$.
Define the set $\L_0 =\{\l\in \L: \l=E_n^+\in \L_1 $ and $\l=E_n^+\in \L_2, \g_n=\es, n\ge 1\}\ss \L$.
In fact with each $\g_n=\es$ we associate two points  $E_n^+\in \L_1$ and $E_n^+\in \L_2$ from the set $\L_0$.
If each gap of $H_0$ is empty, then $\L_0$ is a union of two sets $\{(\pi n)^2, n\in \N\}\ss \L_1$ and $\{(\pi n)^2, n\in \N\}\ss \L_2$.
If each gap of $H_0$ is not empty, then $\L_0=\es$.

\no {\bf Definition of states.} {\it Each zero of $\Xi$ in $\L\sm \L_0$ is a state  of $H$.

\no 1) A state  $\l\in \L_2$  is a resonance.

\no 2) A state   $\l=E_n^\pm , n\ge 0$ is a virtual state.

\no 3) A resonance $\l\in \cup \g_n^{(2)}\ss\L_2$ is an anti-bound state.}

It is known that the gaps $\g_n=\es$ do not give contribution to the states.
Recall that S-matrix for  $H,H_0$ is meromorphic on  $\L$, but it is analytic at the points from $\L_0$ (see \cite{F1}). Roughly speaking there is no difference between the points from $\L_0$ and other points inside the spectrum of $H_0$.

The gap $\g_n^{(1)}\ss\L_1$ is so-called physical gap and the gap $\g_n^{(2)}\ss\L_2$ is so-called non-physical gap.
If $q_0=\int_\R q(x)dx\ne 0$, then $H$ has  precisely one bound state on each physical gap $\ne \es$ and odd number $\ge 1$ of resonances on each non-physical gap $\ne \es$ for $n$ large enough \cite{F1}.

 We explain roughly why antibound states are important. Consider the operator
$H_\t=H_0+\t q$, where $\t\in \R$ is  the coupling constant. If $\t=0$, then $H=H_0$ and there are no eigenvalues, complex resonances and  $H$ has only virtual states, which coincide with each $E_n^\pm , , \g_n\ne \es$, since $\X=2i\O$ on $\L$ at $q=0$ and $\O(E_n^\pm)=0$.  If $\t$ is increasing, then there are eigenvalues, antibound states (close to the end of gaps)  and complex resonances. If $\t$ is increasing again, then there are no new eigenvalues but some complex resonances ($\z\in \C_+\ss\L_2$ and $\ol \z\in \C_-\ss\L_2$) reach some non physical gap and transform into  new antibound states. If $\t$ is increasing again, then
some new antibound states will be virtual states, and then later they will be bound states. Thus if $\t$ runs through $\R_+$, then there is a following transformation: \ resonances $\to $ antibound states $\to $ virtual states  $\to $ bound states $\to $ virtual states...

For each $ n\ge 1$ there exists an unique point $E_n\in [E_n^-,E_n^+]$ such
\[
\lb{dzn}
(-1)^n\D(\sqrt E_n)=\max_{z\in [E_n^-,E_n^+]} |\D(\sqrt\l)|=\cosh h_n, \qq \ {\rm for \ some} \ h_n\ge 0.
\]
Note that if $\g_n=\es$, then  $E_n=E_n^\pm$ and  if  $\g_n\ne \es$, then $E_n\in \g_n$.

\begin{theorem}
\lb{T1}
Let potentials $(p,q)\in L^2(0,1)\ts \cQ_t, t>0$. Then
 $\Xi$ satisfies \er{T1-1} and

\no i) There exist even number $\ge 0$ of states (counted with multiplicity)  on each set $\g_n^{(0)}\ne \es,n\ge 1$, where $\g_n^{(0)}$ is a union of the physical $\ol \g_n^{(1)}\ss \L_1$ and non-physical gap $\ol \g_n^{(2)}\ss\L_2$ (here $\g_n\ne \es$).

\no ii) Let $\|q\|_t=\!\int_0^t\!\!|q(x)|dx$. There are no states in the "forbidden"  domain $\mD_F\sm \cup \ol \g_n^{(2)}\ss\L_2$, where
\begin{multline}
\lb{T1-3}
\mD_F=\{\l\in \L_2: |\l|^{1\/2}>\max \{180e^{2\|p\|_1},C_Fe^{2t|\Im \sqrt \l|}\}\},\ \ C_F=12\|q\|_te^{\|p\|_1+\|q\|_t+2\|p\|_t}.
\end{multline}

\no iii)  In each  $\g_n^{(0)}\ne \es, n\ge 1+{e^{t\pi /2}\/\pi}C_F$ there exists exactly two   simple real states $\l_n^\pm\in \g_n^{(0)}$ such that $E_n^-\le\l_n^-<E_n<\l_n^+\le E_n^+$.

\no iv) Let $\l\in \g_n^{(1)}$ be an eigenvalue of $H$ for some $n\ge 0$, i.e.,
$\Xi(\l)=0$ and let $\l^{(2)}\in \g_n^{(2)}\ss \L_2$ be the same number but on the second sheet $\L_2$.
Then $\l^{(2)}$ is not an anti-bound state, i.e., $\Xi(\l^{(2)})\ne 0$.

\end{theorem}

\no {\bf Remark.} 1) The forbidden domain $\mD_F\cap \C_-$  is similar to the case $p=0$, see \cite{K2}.

2) In Theorem \ref{T2} we show that $\l_n^\pm\to E_n^\pm$ as $n\to \iy$.

Let $\m_n^2, n\ge 1$ be the Dirichlet  spectrum of the equation
$-y''+py=\m^2y$ on the interval $[0,1]$ with the boundary condition $y(0)=y(1)=0$. It is well known that each $\m_n^2\in [E^-_n,E^+_n ], n\ge 1$. Define the coefficients $q_0=\int_\R q(x)dx$ and
\[
\lb{fco}
\wh q_n=\wh q_{cn}+i\wh q_{sn},\qq
\wh q_{cn}=\int_\R q(x)\cos 2\pi nxdx,\qqq
\wh q_{sn}=\int_\R q(x)\sin 2\pi nxdx, \ n\ge 1.
\]
In order to formulate Theorem \ref{T2} we define the functions $\P_{cn}, \P_{sn},
c_n=\cos \f_n, s_n=\sin \f_n, n\ge 1$ (depending from $p\in L^2(0,1)$)  by
\[
\lb{ip1}
 \P_{cn}={E_n^-+E_n^+\/2}-\m_n^2={|\g_n|\/2}c_n,\qq
\P_{sn}={|\g_n|\/2}s_n,\
\qq
\sign s_n=\sign |\vp'(1,\m_n)|.
\]
The identity ${E_n^-+E_n^+\/2}-\m_n^2={|\g_n|\/2}c_n$ defines $c_n=\cos \f_n\in [-1,1]$ and the identity $\sign s_n=\sign |\vp'(1,\m_n)|$ defines $s_n=\sin \f_n$ and  the angle $\f_n\in [0,2\pi)$. Recall the results from \cite{K5}:

{\it The mapping $\P: \cH\to \ell^2\os \ell^2$ given by $\P=((\P_{cn})_1^\iy,(\P_{sn})_1^\iy)$ is a real analytic isomorphism between real Hilbert spaces $\cH=\{p\in L^2(0,1): \int_0^1p(x)dx=0\}$ and $\ell^2\os \ell^2$.}

Let $\#(H,r, X)$ be the total number of state of $H$ in the set
$X\subseteq \L$ having modulus $\le r^2$, each state being counted according to its multiplicity.

\begin{theorem}
\lb{T2}
Let  $(p,q)\in L^2(0,1)\ts \cQ_t$ and let $q_0=\int_\R q(x)dx$. Then

\no  i) Let $\g_n\ne \es$ for some $n$ large enough. Then   $H$ has exactly two  simple states $\l_n^-,\l_n^+\in \g_n^{(0)}$. Moreover, if $\l$ is one of $\l_n^-,\l_n^+$ and satisfies
$(-1)^{n+1}J(\l)>0$ (or $(-1)^{n+1}J(\l)<0$ or $J(\l)=0$), then $\l$ is a bound state (or an anti-bound state or a virtual state).

\no  ii) The following asymptotics  hold true:
\begin{multline}
\lb{T2-1}
\sqrt{\l_n^\pm}=\sqrt{E_n^\pm}\mp {2|\g_n|\/(4\pi n)^3}(\mp q_0-c_n\wh q_{cn}+s_n\wh q_{sn}+ O(1/n))^2,\\
\qqq (-1)^{n+1}J(\l_n^\pm)={|\g_n|\/(2\pi n)^2}(\mp q_0-c_n\wh q_{cn}+s_n\wh q_{sn}+ O(1/n)),
\\
if \qq q_0>0 \ \Rightarrow \qq \l_n^-\in \L_1\ is \ bound  \ state, \qqq
\l_n^+\in \L_2\ is \ anti \ bound  \ state, \\
 if \qq q_0<0 \ \Rightarrow \ \qq
\l_n^-\in \L_2 \ is \ anti \ bound  \ state, \qqq  \l_n^+\in \L_1\ is \ bound  \ state.
\end{multline}

\no iii) Let $\wh q_n=\wh q_{cn}+i\wh q_{sn}=|\wh q_n|e^{i\t_n}, n\ge 1$ and $q_0=0$.
Assume that $|\cos (\f_n+\t_n)|>\ve>0$ and $|\wh q_n|>n^{-\a}$ for  $n$ large enough and for some  $\ve, \a\in (0,1)$. Then the operator $H$ has
$\s_n=1-\sign \cos (\f_n+\t_n)$ bound states in the physical gap $\g_n^{(1)}\ne \es$
and $2-\s_n$ resonances inside the  nonphysical gap $\g_n^{(2)}\ne \es$  for $n$ large enough.

\no iv) For some integer $N_S\in \Z$ the following asymptotics hold true:
\[
\lb{T2-2}
\#(H,r,\L_2\sm\cup\g_n^{(2)})=r{2t+o(1)\/\pi}\qq as \qq r\to \iy,
\]
\[
\lb{T2-3}
\#(H,r,\cup \G_n)=\#(H_0,r,\cup \G_n)+2N_S \qq as \qq r\to \iy, \qq r^2\notin \cup\ol\g_n.
\]
\end{theorem}

\no{\bf Remark}. 1) First term in the asymptotics \er{T2-2} does not depend
on the periodic potential $p$.  Recall that asymptotics \er{T2-2} for the case $p=0$ was obtained
by Zworski \cite{Z}.

\no 2) The main difference between the distribution of the resonances for the case $p\ne \const$
and $p=\const$ is the bound states and anti bound states in high  energy gaps, see  (iii) and \er{T2-1}.

\no 3) Assume that a potential $u\in L^2(\R)$ is compactly supported, $\supp u\ss (0,t)$  and satisfies $|\wh u_{n}|=o(n^{-\a})$ as $n\to \iy$.
Then in the case (iii) the operator $H+u$ has also $1-\sign \cos (\f_n+\t_n)$ bound states in each gap $g_n\ne \es$ for $n$ large enough.

\no 4)  In the proof of \er{T2-2} we use the Paley Wiener type Theorem from \cite{Fr}, the Levinson Theorem
(see Sect. 3) and a priori estimates from \cite{KK}, \cite{M}.

\no 5) Assume that  $\m_n^2=E_n^-$ or $\m_n^2=E_n^+$ with $|\g_n|>0$ for all
$n\in \N_0$, where $ \N_0\ss \N$ is some infinite subset. Then
for all large $n\in \N_0$ the following asymptotics hold true:
\[
\lb{T4-1}
\sqrt{\l_n^\pm}=\sqrt{E_n^\pm}\mp {2|\g_n|\/(4\pi n)^2}\rt(\mp q_0-c_n\wh q_{cn}+ {O(1)\/n}\rt)^2
\qqq \as \qq n\to \iy.
\]
 There is no  $\wh q_{sn}$ in asymptotics \er{T4-1}.
Let in addition $q_0=0$ and $|\wh q_{cn}|>n^{-\a}$ for all $n\in \N_0$ and some $\a\in (0,1)$. Then
$H$ has exactly $\s_n=1-\sign c_n\wh q_{cn}$ eigenvalues in each gap $\g_n^{(1)}$
and  $2-\s_n$ resonances inside each gap $\g_n^{(2)}\ne \es$ for $n\in \N_0$ large enough.

\no 6) Recall that if $p$ is even i.e., $p\in L_{even}^2(0,1)=\{p\in L^2(0,1),
p(x)=p(1-x), x\in (0,1)\}$, iff $s_n=0$ ( or  $\m_n^2\in \{E_n^-,E_n^+\}$,  or
$c_n\in \{\pm 1\}$) for all $n\ge 1$, see \cite{GT}, \cite{KK1}.

Consider  some inverse problems for the operator $H$.

\begin{theorem}
\lb{T3}
i)  Let the operator $H_0$ have infinitely many gaps $\g_n\ne \es$
for some $p\in L^2(0,1)$. Then for any sequence $\s=(\s_n)_{1}^\iy$,
where $\s_n\in \{0,2\}$,  there exists some potential $q\in \cQ^t,t>0$ such that $H$ has exactly $\s_n$ eigenvalues in each gap $\g_n^{(1)}\ne \es$  and  $2-\s_n$ resonances inside each gap $\g_n^{(2)}\ne \es$ for $n$ large enough.

\no ii)  Let  $q\in \cQ_t,t>0$ satisfy   $q_0=0$ and let $|\wh q_{cn}|>n^{-\a}$ for all $n$ large enough and some $\a\in (0,1)$. Let $\s=(\s_n)_{1}^\iy$, where $\s_n\in \{0,2\}$ and
infinitely many $\d_n>0$. Then for any sequence  $\d=(\d_n)_1^\iy\in \ell^2$ where each $\d_n\ge 0, n\ge 1$
there exists a potential $p\in L^2(0,1)$  such that each gap length  $|\g_n|=\d_n, n\ge 1$.
Moreover, $H$ has exactly $\s_n$ eigenvalues in each physical gap $\g_n^{(1)}\ne \es$
and   $2-\s_n$ resonances inside each non-physical gap $\g_n^{(2)}\ne \es$ for $n$ large enough.
\end{theorem}

{\bf Remark .} 1)  In the proof we use results from the inverse spectral theory from \cite{K5}.

\no 2)  We have $\a\in ({1\/2},1)$, if all gaps are open.

A lot of papers are devoted to the resonances for the Schr\"odinger
operator with $p=0$, see \cite{Fr}, \cite{H}, \cite{K1}, \cite{K2}, \cite{S}, \cite{Z}  and
references therein. Although resonances have been studied in many settings, but there are relatively
few cases where the asymptotics of the resonance counting function are known, mainly one dimensional case  \cite{Fr}, \cite{K1}, \cite{K2}, \cite{S}, and \cite{Z}. We recall that
Zworski [Z] obtained the first results about the distribution of
resonances for the Schr\"odinger operator with compactly supported
potentials on the real line. The author obtained the characterization (plus uniqueness and
recovering) of $S$-matrix for the Schr\"odinger operator with a compactly supported potential on the
real line \cite{K2} and the half-line \cite{K1}, see also \cite{Z1}, \cite{BKW} about uniqueness.

For the Schr\"odinger operator on the half line the author \cite{K3}  obtained the stability results:

(i) If $\vk^0=(\vk^0)_1^\iy$ is a sequence of zeros (eigenvalues
and resonances) of the Jost function for some real compactly
supported potential $q_0$ and $\vk-\vk^0\in\ell_\ve^2$ for some
$\ve>1$, then $\vk$ is the sequence of zeroes of the Jost function
for some unique real compactly supported potential.

(ii) The measure associated with the zeros of the Jost
function is the Carleson measure, the sum $\sum
(1+|\vk_n^0|)^{-\a}, \a>1$ is estimated in terms of the $L^1$-norm of the potential $q^0$.

Brown and Weikard \cite{BW} considered  the Schr\"odinger operator $-y''+(p_A+q)y$ on the half-line,
where $p_A$ is an algebro-geometric potentials
and $q$ is a compactly supported potential. They proved that the zeros of the Jost function determine $q$ uniquely.

Christiansen \cite{Ch} considered resonances associated to the Schr\"odinger
operator $-y''+(p_{S}+q)y$ on the real line, where $p_S$ is a step potential. She determined asymptotics of the resonance-counting function. Moreover, she obtained  that the
resonances determine $q$ uniquely.

In the proof of theorems we use properties of the quasimomentum from \cite{KK}, \cite{M},
a priori estimates from \cite{KK}, \cite{M},   and
results from the inverse theory for the Hill operator from \cite{K5}, see \er{mz}.

The plan of the paper is as follows. In Section 2 we describe the preliminary results about fundamental solutions. In Sections 3 we study the function $\x$ and prove Theorem \ref{T1}. In Sections 4 we prove the main
Theorem \ref{T2}-\ref{T3}, which are crucial for
Section 3 and 4.

\section {Preliminaries }
\setcounter{equation}{0}

We will work with the momentum $z=\sqrt \l$, where $\l$ is an energy and recall
that $\sqrt 1=1$.  If $\l\in \g_n,  n\ge 1$, then $z\in g_{\pm n}$ and if $\l\in \g_0=(-\iy,0)$, then  $z\in g_0^\pm=i\R_\pm$, where the momentum gaps $g_n$ are given by
\[
\lb{2}
g_n=(e_n^-,e_n^+)=-g_{-n},\qq e_n^\pm=\sqrt{E_n^\pm}>0,\qq
n\ge 1,\qq  and  \qq \D(e_{n}^{\pm})=(-1)^n.
\]
Introduce the cut domain (see Fig.2)
$$
\cZ_0=\C\sm \cup \ol g_n.
$$
\begin{figure}
\tiny
\unitlength=1mm
\special{em:linewidth 0.4pt}
\linethickness{0.4pt}
\begin{picture}(120.67,34.33)
%coordinate lines
\put(20.33,21.33){\line(1,0){100.33}}
\put(70.33,10.00){\line(0,1){24.33}}
\put(69.00,19.00){\makebox(0,0)[cc]{$0$}}
\put(120.33,19.00){\makebox(0,0)[cc]{$\Re z$}}
\put(67.00,33.67){\makebox(0,0)[cc]{$\Im z$}}
%slits
\put(81.33,21.33){\linethickness{2.0pt}\line(1,0){9.67}}
\put(100.33,21.33){\linethickness{2.0pt}\line(1,0){4.67}}
\put(116.67,21.33){\linethickness{2.0pt}\line(1,0){2.67}}
\put(60.00,21.33){\linethickness{2.0pt}\line(-1,0){9.33}}
\put(40.00,21.33){\linethickness{2.0pt}\line(-1,0){4.67}}
\put(24.33,21.33){\linethickness{2.0pt}\line(-1,0){2.33}}
%endpoints
\put(81.67,24.00){\makebox(0,0)[cc]{$e_1^-$}}
\put(91.00,24.00){\makebox(0,0)[cc]{$e_1^+$}}
\put(100.33,24.00){\makebox(0,0)[cc]{$e_2^-$}}
\put(105.00,24.00){\makebox(0,0)[cc]{$e_2^+$}}
\put(115.33,24.00){\makebox(0,0)[cc]{$e_3^-$}}
\put(120.00,24.00){\makebox(0,0)[cc]{$e_3^+$}}
\put(59.33,24.00){\makebox(0,0)[cc]{$-e_1^-$}}
\put(50.67,24.00){\makebox(0,0)[cc]{$-e_1^+$}}
\put(40.33,24.00){\makebox(0,0)[cc]{$-e_2^-$}}
\put(34.67,24.00){\makebox(0,0)[cc]{$-e_2^+$}}
\put(26.00,24.00){\makebox(0,0)[cc]{$-e_3^-$}}
\put(19.50,24.00){\makebox(0,0)[cc]{$-e_3^+$}}
\end{picture}
\caption{The cut domain $\cZ=\C\sm \cup g_n$ and the slits $g_n=(e_n^-,e_n^-)$ in the $z$-plane.}
\lb{z}
\end{figure}

Below we will use the momentum
variable $z=\sqrt \l$ and the corresponding the  Riemann surface $\cM$, which is more convenient for us, than the Riemann surface $\L$.
Slitting the n-th momentum gap  $g_n$ (suppose it is nontrivial) we obtain a cut $G_n$ with an upper   $g_n^+$ and lower rim  $g_n^-$. Below we will identify this cut $G_n$ and the union of of the upper rim  (gap) $\ol g_{n}^+$ and the lower rim (gap) $\ol g_{n}^{\ -}$, i.e.,
\[
G_{n}=\ol g_{n}^+\cup \ol g_{n}^-,\ \ where \ \ g_{n}^\pm =g_n\pm i0.
%\qq If \ z\in g_n \Rightarrow z\pm i0\in g_n^\pm.
\]
Also we will write  $z\pm i0\in g_n^\pm$, if $z\in g_n$.

In order to construct the  Riemann surface $\cM$
we take the cut domain $\cZ_0=\C\sm \cup \ol g_n$ and identify (i.e. we glue) the upper rim $g_{n}^+$ of the slit $G_n$ with the upper rim $g_{-n}^+$  of the slit $G_{-n}$ and correspondingly the lower rim $g_{n}^-$ of the slit $G_n$ with the  lower  rim $g_{-n}^-$  of the slit $G_{-n}$ for all nontrivial gaps.
%We have the one-to-one mapping from $\L$
%onto $\cM$ given by $z=\sqrt E$.
%Such surface $\cM$ was introduce by Firsova \cite{F3}.

The mapping $z=\sqrt \l$ from $\L$ onto $\cM$ is one-to-one and onto.
The gap $\g_n^{(1)}\ss \L_1$ maps onto $g_n^+\ss \cM_1$ and the gap $\g_n^{(2)}\ss \L_2$ maps onto $g_n^-\ss \cM_2$.
From a physical point of view, the upper rim  $g_{n}^+$ is a physical gap
and the lower rim  $g_{n}^-$ is a non physical gap.
Moreover, $\cM\cap\C _+=\cZ_0\cap\C _+$ plus all physical gaps $g_{n}^+$ is a so-called physical "sheet" $\cM_1$
and $\cM\cap\C _-=\cZ_0\cap\C _-$ plus all non physical gaps $g_{n}^-$
is a so-called non physical "sheet" $\cM_2$. The set (the spectrum) $\R\sm \cup g_n$ joints the first and second sheets.

We introduce the quasimomentum $k(\cdot )$ for $H_0$ by
$k(z)=\arccos \D(z),\ z \in \cZ_0=\C\sm \cup \ol g_n$. The function
$k(z)$ is analytic in $z\in\cZ_0$ and satisfies:
\[
\lb{pk}
(i)\qq k(z)=z+O(1/z)\qq  as \ \ |z|\to \iy, \qq and \qq (ii)\qq  \Re k(z\pm i0) |_{[e_n^-,e_n^+]}=\pi n,\qq \ n\in \Z,
\]
and $\pm \Im k(z)>0$ for any $z\in \C_\pm$, see \cite{M}, \cite{KK}. The function $k(\cdot)$ is analytic
on $\cM$ and satisfies $\sin k(z)=\O(z^2), z\in \cM$. Moreover, the quasimomentum $k(\cdot)$ is a conformal mapping from $\cZ_0$ onto the quasimomentum domain $\cK=\C\sm \cup \G_n$, where the slit $\G_n=(\pi n-ih_n,\pi n+ih_n)$.
The height $h_n$ is defined by the equation $\cosh h_n=(-1)^n\D(e_n)$, where $e_n\in [e_n^-,e_n^+]$
and $\D'(e_n)=0$.
The Floquet solutions $\p_{\pm}(x,z), z \in \cZ_0$ of $H_0$ is given by
\[
\lb{3}
\p_\pm(x,z)=\vt(x,z)+m_\pm(z)\vp(x,z),\ \
m_\pm={\b\pm i\sin k\/ \vp(1,\cdot)},\ \
\b={\vp'(1,\cdot)-\vt(1,\cdot)\/2},
\]
where $\vp(1,z)\p^+(\cdot,z)\in L^2(\R_+)$ for all $z\in\C_+\cup\cup_{} g_n$.
If $p=0$, then $k=z$ and $\p_\pm(x,z)=e^{\pm izx}$.

\begin{figure}
\tiny
\unitlength=1mm
\special{em:linewidth 0.4pt}
\linethickness{0.4pt}
\begin{picture}(120.67,34.33)
%coordinate lines
\put(20.33,20.00){\line(1,0){102.33}}
\put(71.00,7.00){\line(0,1){27.00}}
\put(70.00,18.67){\makebox(0,0)[cc]{$0$}}
\put(124.00,18.00){\makebox(0,0)[cc]{$\Re k$}}
\put(67.00,33.67){\makebox(0,0)[cc]{$\Im k$}}
%slits
\put(87.00,15.00){\linethickness{2.0pt}\line(0,1){10.}}
\put(103.00,17.00){\linethickness{2.0pt}\line(0,1){6.}}
\put(119.00,18.00){\linethickness{2.0pt}\line(0,1){4.}}
\put(56.00,15.00){\linethickness{2.0pt}\line(0,1){10.}}
\put(39.00,17.00){\linethickness{2.0pt}\line(0,1){6.}}
\put(23.00,18.00){\linethickness{2.0pt}\line(0,1){4.}}
\put(85.50,18.50){\makebox(0,0)[cc]{$\pi$}}
\put(54.00,18.50){\makebox(0,0)[cc]{$-\pi$}}
\put(101.00,18.50){\makebox(0,0)[cc]{$2\pi$}}
\put(36.00,18.50){\makebox(0,0)[cc]{$-2\pi$}}
\put(117.00,18.50){\makebox(0,0)[cc]{$3\pi$}}
\put(20.00,18.50){\makebox(0,0)[cc]{$-3\pi$}}
%endpoints
\put(87.00,26.00){\makebox(0,0)[cc]{$\pi+ih_1$}}
\put(56.00,26.00){\makebox(0,0)[cc]{$-\pi+ih_1$}}
\put(103.00,24.00){\makebox(0,0)[cc]{$2\pi+ih_2$}}
\put(39.00,24.00){\makebox(0,0)[cc]{$-2\pi+ih_2$}}
\put(119.00,23.00){\makebox(0,0)[cc]{$3\pi+ih_3$}}
\put(23.00,23.00){\makebox(0,0)[cc]{$-3\pi+ih_3$}}
\end{picture}
\caption{The domain $\cK=\C\sm \cup \G_n$, where the slit $\G_n=(\pi n-ih_n,\pi n+ih_n)$}
\lb{k}
\end{figure}

The function $\sin k$ and each function $\vp(1,\cdot)\p_{\pm}(x,\cdot), x\in \R$ are analytic on the Riemann surface $\cM$.
%Moreover, $\sin k(z) |q|^{1\/2}(H_0-z^2)^{-1}|q|^{1\/2}$ is analytic %operator-valued function on $\cM$ (see Sect.2) and belongs to the trace class.
%Note that the functions $\D(\sqrt E)$ and $\vt(x,\sqrt E), \vp(x,\sqrt E),.. $ %are entire in $E\in \C$.
%The functions $\sin k(\sqrt E)$ and the functions $\vp(1,z)\p_{\pm}(x,\sqrt E)$ %are analytic on the two-sheeted
%Riemann surface $\L$ of $(\D^2(\sqrt E)-1)^{1\/2}$ obtained by joining the upper %and lower rims of two copies
%of the cut plane $\C\sm\s_{ac}(H_0)$ in the usual (crosswise) way.
%The first sheet $\L_1$ is a so-called physical sheet, where
%the resolvent $(H_0-E)^{-1}$ is analytic operator-valued function
%and all eigenvalues of $H$ belong to $\L_1$.
%The operator $(H_0-E)^{-1}$ is not bounded on the second sheet
%and all resonances belong to $\L_2$.
Recall that the Floquet solutions $\p_\pm(x,z), (x,z)\in \R\ts \cM$ satisfy (see \cite{T})
\[
\lb{f1}
\p_\pm(0,z)=1, \qq \p_\pm(0,z)'=m_\pm(z),
\qq \p_\pm(1,z)=e^{\pm ik(z)}, \qq  \p_\pm(1,z)'=e^{\pm ik(z)}m_\pm(z),
\]
\[
\lb{f2}
\p_\pm(x,z)=e^{\pm ik(z)x}(1+O(1/z)) \qq as \qq |z|\to \iy, \qq (x,z)\in \R\ts \cZ_\ve,
\]
where the set $\cZ_\ve =\{z\in \cZ_0: \dist \{z,g_n\}>\ve, g_n\ne \es,  n\in \Z\},\ \ve>0$.
Below we need the simple identities
\[
\lb{LD0}
\b^2+1-\D^2=1-\vp'(1,\cdot)\vt(1,\cdot)= -\vp(1,\cdot)\vt'(1,\cdot).
\]
Let $\cD_r(z_0)=\{|z-z_0|<r\}$ be a disk for some $r>0$.
Recall the well known properties of the function $m_\pm$ (see \cite{T}).

\begin{lemma}
\lb{Tm}
i) $\Im m_+ (z)>0$ for all $(z,n)\in (z_{n-1}^+,z_{n}^-)\ts \N$ and the following asymptotics hold true:
\[
\lb{Tm-1}
m_\pm (z)=\pm iz+O(1) \qq as \qq |z|\to \iy, \qq z\in \cZ_\ve,\ve >0.
\]

\no ii) Let $g_n=\es$ for some $n\in \Z$. Then
the functions $\sin k(\cdot), m_\pm$ are analytic in the disk
$B(\m_n,\ve)\ss\cZ_0$ for some $\ve>0$ and the functions $\sin k(z)$
and $\vp(1,z)$ have the simple zero at $\m_n$. Moreover, $m_\pm$ satisfies
\[
\lb{Tm-2}
m_\pm (\m_n)={\b'(\m_n)\pm i(-1)^nk'(\m_n)\/\pa_z\vp(1,\m_n)},
\qq \Im m_\pm (\m_n)\ne 0.
\]

\no iii) If the function $m_+$ has a pole at $\m_n+i0$
 for some $n\in \N$, then  $k(\m_n+i0)=\pi n+ih_{sn}$ and
\begin{multline}
\lb{Tm-31}
\qq h_{sn}>0,\qq
\b(\m_n)=i\sin k(\m_n+i0)= -(-1)^n\sinh h_{sn},\qq m_+\in \mA(\m_n-i0),\\
 \qq m_+(\m_n+z)={\r_n+O(z)\/z} \qq as \ z\to 0, \ z\in \C_+,\
\qq \r_n={-2\sinh |h_{sn}|\/(-1)^n\pa_z\vp(1,\m_n)}<0.
\end{multline}

\no iv)  If the function $m_+$ has a pole at $\m_n-i0$
for some $n\in \N$, then $k(\m_n-i0)=\pi n+ih_{sn}$ and
\begin{multline}
\lb{Tm-32}
h_{sn}<0,\qq
\b(\m_n)=-i\sin k(\m_n-i0)=(-1)^n\sinh h_{sn},\ m_+\in \mA(\m_n^+),\\
 \qq m_+(\m_n+z)={\r_n+O(z)\/z}\qq as \ z\to 0, \ z\in \C_-.
\end{multline}

\no iv)  $\m_n=e_n^-$ or $\m_n=e_n^+$ (here $e_n^-\ne e_n^+$) for some $n\ne 0$ iff
\[
m_+(\m_n+z)={\r_n^\pm+O(z)\/\sqrt z}\qq as \ z\to 0, \ z\in \C_+,\qq\ {\rm some} \ \const \
\r_n^\pm\ne 0.
\]
\end{lemma}

 The following asymptotics hold true as $n\to \iy$ (see \cite{PT}, \cite{K5}):
\[
\lb{sde}
\m_n=\pi n+\ve_n(p_{c0}-p_{cn}+O(\ve_n)),
\qqq
\ve_n={1/2\pi n},
\]
\[
\lb{ape}
e_n^\pm=\pi n+\ve_n(p_0\pm |p_n|+O(\ve_n)), \qq \qq
p_n=\int_0^1p(x)e^{-i2\pi nx}dx=p_{cn}-ip_{sn}.
\]
The equation $-f''+(p+q)f=z^2 f,\ z\ne e_n^\pm$
has unique solutions $f_\pm (x,z)$ such that
\[
\lb{bcf}
f_+(x,z)=\p_+(x,z), \ x\ge t, \qqq {\rm and} \ \ \ \
f_-(x,z)=\p_-(x,z), \ x\le 0,
\]
and $f_+(x,z)=\ol f_+(x,-z), z^2\in\s(H_0)\sm \{e_n^\pm, n\in \Z\}$. This yields
\[
f_+(x,z)=b(z)f_-(x,z)+a(z)f_-(x,-z),\qq a={w\/w_0},\qq b={s\/w_0},
\]
where $z^2\in\s(H_0)\sm \{e_n^\pm, n\in \Z\}$ and
\begin{multline}
\ w=\{f_-,f_+\},  \ \ \ \
s=\{f_+(x,z),f_-(x,-z)\},\
w_0=\{\p_-,\p_+\}={2i\sin k\/\vp(1,\cdot)}
\end{multline}
and $\{f, g\}=fg'-f'g$ is the Wronskian.
The scattering matrix $\cS_M$ for $H, H_0$ is given by
\[
\cS_M (z)\ev \ma a(z)^{-1}& r_-(z)\\ r_+(z)&a(z)^{-1}\am ,
\ \ \ \  r_{\pm}={s(\mp z)\/ w(z)}=\mp {b(\mp z)\/ a(z)}, \ \ \ \  z^2\in \s_{ac}(H),
\lb{1.7}
\]
where $1/a$ is the transmission coefficient and $r_{\pm}$ is the reflection coefficient. We have the following identities from \cite{F1}, \cite{F3}:
\[
\lb{iab}
|a(z)|^2=1+|b(z)|^2,\qqq z^2\in \s_{ac}(H),
\]
\[
\lb{iab1}
a(z)=D(z^2), \qqq z\in \cZ_0.
\]
The functions $a,b, s,w $ are analytic in $\cZ_0$ and real on $i\R$. Then the following identities hold true:
\[
a(-z)=\ol a(\ol z), \qq w(-z)=\ol w(\ol z), \qq s(-z)=\ol s(\ol z),\qq w_0(-z)=\ol w_0(\ol z),\qq  z\in \cZ_0.
\]

Let $\vp(x,z,\t), \ (z,\t)\in \C\ts \R$   be the solutions of the equation
\[
\lb{x+t}
-\vp''+p(x+\t)\vp=z^2 \vp, \qq \ \vp(0,z,\t)=0,\qq \vp'(0,z,\t)=1.
\]
The function $\vp(1,z,x)$ for all $(x,z)\in R\ts \C$ satisfies the following identity (see \cite{Tr})
\[
\lb{if}
\vp(1,\cdot,\cdot)=\vp(1,\cdot)\vt^2-\vt'(1,\cdot)\vp^2+2\b\vp\vt=
\vp(1,\cdot)\p_-\p_+.
\]
Let  $\wt\vt, \wt\vp$ be the solutions of the equations
$-y''+(p+q)y=z^2y, z\in \C$ and satisfying
\[
\lb{wtc}
\wt\vp(x,z)=\vp(x,z),\qqq \wt\vt(x,z)=\vt(x,z) \ \ for \ \ all \ x\ge t.
\]
A solution of the equation $-y''+(p-z^2)y=f, y(0)=y'(0)=0$  has the form
$
y=\int_0^x\vp(x-\t,z,\t)f(\t)d\t.
$
Hence the solutions $\wt\vt, \wt\vp$ and $f_+$  of the equation
$-y''+(p+q)y=z^2y$ satisfy the equation
\[
\lb{ep}
y(x,z)=y_0(x,z)-\int_x^t\vp(x-\t,z,\t)q(\t)y(\t,z)d\t, \qqq x\le t,
\]
where $y$ is one from  $\wt\vt, \wt\vp$ and $f_+$; $y_0$ is the corresponding function from  $\vt, \vp$ and $\p_+$.
For each $x\in \R$ the functions  $\wt\vt (x,z), \wt\vp (x,z)$ and
$\vt (x,z), \vp (x,z)$ are entire in $z\in\C$ and satisfy
\begin{multline}
\lb{efs}
\max \{||z|_1\wt\vp(x,z)|, \ |\wt\vp'(x,z)| , |\wt\vt(x,z)|,
{1\/|z|_1}|\wt\vt'(x,z)|    \} \le X_1=e^{|\Im z||2t-x|+\|q\|_t+\|p\|_t+\int_x^t|p(\t)|d\t},\\
|z|_1=\max\{1, |z|\},\qqq
|\wt\vt(x,z)-\vt(x,z)|\le {X_1\/|z|}\|q\|_t,\qq
|\wt\vp(x,z)-\vp(x,z)|\le {X_1\/|z|^2}\|q\|_t,
\end{multline}
where $\|p\|_t=\int_0^t|p(s)|ds$   and
\begin{multline}
\lb{efs1}
\max \{|z|_1|\vp(x,z)|, \ |\vp'(x,z)| , |\vt(x,z)|,
{1\/|z|_1}|\vt'(x,z)|    \} \le X=e^{|\Im z|x+\|p\|_x},\\
|\vp(x,z)-{\sin zx\/z}|\le {X\/|z|^2}\|p\|_x,
\qq |\vt(x,z)-{\cos zx}|\le {X\/|z|}\|p\|_x,
\end{multline}
for  $(p,x,z)\in L_{loc}^1(\R)\ts \R\ts \C$, see  [PT]. These estimates
yield
\begin{multline}
\lb{asb} \b(z)=\int_0^1{\sin z(2x-1)\/z}p(x)dx+{O(e^{|\Im
z|})\/z^2},\\ \b'(z)=\int_0^1{\cos z(2x-1)\/z} p(x)(2x-1)
dx+{O(e^{|\Im z|})\/z^2}\qq as \qq |z|\to \iy.
\end{multline}

\begin{lemma}
\lb{T22}
For all $z\in \cZ_0$ the following identities and asymptotics  hold true:
\[
\lb{T22-1}
f_+(\cdot,z)=\wt\vt (\cdot,z)+m_+(z)\wt\vp (\cdot,z),
\]
\begin{multline}
\lb{T22-2}
f_+(0,z)\!\!=1\!+\!\int_0^t\!\vp(x,z)q(x)f_+(x,z)dx,\
f_+'(0,z)\!\!=\! m_+(z)-\!\int_0^t\! \vt(x,z)q(x)f_+(x,z)dx,
\end{multline}
\[
\lb{T22-3} f_+(x,z)=e^{ik(z)x}+e^{(2t-x)\gJ}O(1/z),\qqq \qq
\gJ={|\Im z|- \Im z\/2},    \qq x\in [0,t],
\]
\[
\lb{T22-4}
f_+(0,z)=1+e^{2t\gJ}O(1/z),   \qq f_+'(0,z)=iz+O(1)+e^{2t\gJ}o(1)
\]
as $|z|\to \iy, z\in \cZ_\ve, \ve>0$, where $\wh q(z)=\int_0^tq(x)e^{2izx}dx, \ z\in \C$.

\end{lemma}
\no{\bf Proof.} Using \er{wtc}, \er{bcf} we obtain \er{T22-1}.
Using the identity $\vp(x,\cdot,s)=\vt(t,\cdot)\vp(x+s,\cdot)-\vp(s,,\cdot)\vt(x+s,\cdot)$, we obtain
$\vp(-s,\cdot,s)=-\vp(s,\cdot)$ and $\vp'(-s,\cdot,s)=\vt(s,\cdot)$. Substituting the last
identities into \er{ep} we get \er{T22-2}.

Standard iteration arguments (see \cite{PT}) for the equation \er{ep}  give \er{T22-3}.

Substituting \er{T22-3} into \er{ep}  we obtain \er{T22-4}.
\BBox

\begin{lemma}
\lb{T23}
i) For each $z\in \cZ_0$ the following identities hold true:
\[
\lb{T23-1}
 w(z)=f_+'(0,z)-m_-(z)f_+(0,z)=w_0(z)-\int_0^tq(x)\p_-(x,z)f_+(x,z)dx,
\]
\[
\lb{T23-2}
s(z)=f_+(0,z)m_+(z)-f_+'(0,z)=\int_0^tq(x)\p_+(x,z)f_+(x,z)dx.
\]
 The functions $\x(z)=2i\sin k(z)a(z)$ and $s(\cdot)$ have the following asymptotics:
\[
\lb{T23-3}
\x(z)=2i\sin z(1+O(e^{2t\gJ}/z)) ,\qqq s(z)=O(e^{2t\gJ}),
\]
as $|z|\to \iy, z\in \cZ_\ve, \ve>0$.

\no ii) The function $s(\cdot)$ has exponential type $\r_\pm$ in the half plane $\C_\pm$, where $\r_+=0, \r_-=2t$.

\end{lemma}
\no{\bf Proof.} We have $w=\{f^-,f^+\}=\p^-{f^+}'-m^-f^+|_{x=0}$,
which yields the identity in \er{T23-1}.
Substituting \er{T22-2} into \er{T23-1} we obtain \er{T23-2}.
Asymptotics from Lemma \ref{T22} and \er{Tm-1}  imply \er{T23-3}.

ii) We show $\r_-=2t$.  Due to \er{T23-3}, $s$ has exponential type $\r_-\le 2t$.
The decompositions $f_+=e^{ixz}(1+h)$ and $\p_+=e^{ixz} (1+\e)$ give $(1+h)(1+\e)=1+T, T=h+\e +\e h$ and
\[
\lb{esff}
s(z)=\int_0^tq(x)\p^+(x,z)f^+(x,z)dx=\int_0^tq(x)e^{i2xz}(1+T(x,z))dx,
\qqq  \qq z\in\cZ_\ve.
\]

Asymptotics \er{f2}, \er{efs1}, \er{T22-3} and $k(z)=z+O(1/z)$ as $|z|\to \iy$ (see \cite{KK}) yield
\[
\lb{esff1}
\e(x,z)=O(1/z),\qq \qqq h(x,z)=e^{2(t-x)|\Im z|}O(1/z)\qq as \ |z|\to \iy,\qq z\in \cZ_\ve.
\]
We need the following variant of the Paley Wiener type Theorem from \cite{Fr}:

\no {\it let $q\in \cQ_t$ and let each $G(x,z), x\in [0,t]$ be analytic for $z\in\C_-$ and $G\in L^2((0,t)dx,\R dz)$.
Then $\int_0^te^{2izx}q(x)(1+G(x,z))dx$ has exponential type at least $2t$ in $\C_-$.}

We can not apply this result to the function  $T(x,z), z\in\C_-$, since $m_+(z)$ may have a singularity
at $\m_n-i0\in \ol g_n^-$ if  $g_n\ne \es$.  But we can use this result for the function  $T(x,z-i), z\in\C_-$,
since \er{esff}, \er{esff1} imply $\sup _{x\in [0,1]}|T(x,-i+\t)|=O(1/\t)$ as $\t\to \pm\iy$.
Then the function $s(z)$ has exponential type $2t$ in the half plane $\C_-$.
The proof for $\r_+=0$ is similar.
\BBox

\section {Analyze of the function $\x$}
\setcounter{equation}{0}

Below we need the identities and the asymptotics as $n\to \iy$ from \cite{KK}:
\[
\lb{35} (-1)^{n+1}i\sin k(z)=\sinh v(z)=\pm |\D^2(z)-1|^{1\/2}>0\qq
\ all \qq z\in  g_n^\pm,
\]
\[
\lb{pav}
v(z)=\pm |(z-e_n^-)(e_n^+-z)|^{1\/2}(1+O(n^{-2})),\qq
\sinh v(z)=v(z)(1+O(|g_n|^2),\qq  z\in \ol g_n^\pm.
\]
Let $\n_n^2, n\ge 1$ be the Neumann  spectrum of the equation
$-y''+py=\n^2y$ on the interval $[0,1]$ with the boundary condition $y(0)=y(1)=0$. It is well known that each $\n_n^2\in [E^-_n,E^+_n ], n\ge 1$.

\begin{lemma}
\lb{T31}
Let $p\in L^1(0,1)$.

\no i) Then the following asymptotics hold true uniformly for $z\in [e_n^-,e_n^+]$ as $n\to \iy$:
\[
\lb{T31-1}
\vp(1,z)=(-1)^n{(z-\m_n)\/\pi n}(1+O(1/n)),\qq
\]
\[
\lb{T31-2}
-{\vt'(1,z)\/z^2}=(-1)^n{(z-\n_n)\/\pi n}(1+O(1/n)).
\]
 ii) Let $z\in g_n$ and $e_n^-\ge R_p=8e^{\|p\|_1}$.
 Then the following estimates hold true (here $\dot u=\pa_z u$)
\[
\lb{T31-3}  |g_n|^2\le {8e^{\|p\|_1}\/|z_n|}<1,
\]
\[
\lb{T31-4}
|\dot \vp(1,z)|\le {3e^{\|p\|_1}\/2|z|},
\qqq |\vp(1,z)|_{z\in g_n} \le |g_n|{3e^{\|p\|_1}\/2|z|},\qq
\]
\[
\lb{T31-5}
|\dot \vt'(1,z)|\le |z|{3e^{\|p\|_1}\/2},\qqq
|\vt'(1,z)|\le |g_n||z|{3e^{\|p\|_1}\/2},
\]
\[
\lb{T31-6}
|\dot\b(z)|\le  {3e^{\|p\|_1}\/2|z|}, \qqq |\b(e_n^\pm)|\le |g_n|{3e^{\|p\|_1}\/2},\qqq
|\b(z)|\le  |g_n|{9e^{\|p\|_1}\/4|z|}.
\]
\no iii) In each disk $\cD_{\pi\/4}(\pi n)\ss \cD=\{|z|>32e^{2\|p\|_1}\}$
there exists exactly one momentum gap $g_n$ of the operator $H_0$. Moreover, if
$g_n,g_{n+1}\ss\cD$, then $e_{n+1}^--e_{n}^+\ge \pi$.
\end{lemma}
\no {\bf Proof.} i) We have the Taylor series $\vp(1,z)=\dot\vp(1,\m_n)\t+\ddot\vp(1,\m_n+\a \t){\t^2\/2}$
for any $z\in [e_n^-, e_n^+]$ and some $\a\in [0,1]$, where $\t=z-\m_n$ and $\dot \vp=\pa_z \vp$.
Asymptotics \er{efs1} give $\dot\vp(1,\m_n)=2(-1)^n(1+O(1/n))/(2\pi n)$
and $\ddot\vp(1,\m_n+\a \t)\t=O(n^{-2})$, which yields \er{T31-1}.
Similar arguments imply \er{T31-2}.

ii) Using $|\D(z)-\cos z|\le {e^{\|p\|_1}\/|z|}$, for all $|z|\ge 2$, we obtain
$$
{h_n^2\/2}\le \cosh h_n-1=|\D(z_n)|-1\le {e^{\|p\|_1}\/|z_n|}.
$$
Then the estimate $|g_n|\le 2h_n$ (see \cite{KK}) gives  \er{T31-3}.

Due to \er{efs1}, the function $f=z\vp(1,z)$ has the estimate $|f(z)|\le C_0=e^{\|p\|_1}, z\in \R$. Then
the Bernstein inequality gives $|\dot f(z)|=|\vp(1,z)+z\dot \vp(1,z)|\le C_0, z\in \R$, which yields
$|\dot \vp(1,z)|\le {C_0\/|z|}(1+{1\/|z|}), z\in \R$. Moreover, we obtain $|\vp(1,z)|\le |g_n|\max_{z\in g_n} |\dot \vp(1,z)|\le |g_n|{3C_0\/2|z|}$.

The proof of \er{T31-5} and the estimate $|\dot\b(z)|\le  {3e^{\|p\|_1}\/2|z|}$ is similar. Identity \er{LD0} gives

$\b^2(e_n^\pm)=-\vp(1,e_n^\pm)\vt'(1,e_n^\pm)$. Then \er{T31-4}, \er{T31-5} imply $|\b(e_n^\pm)|\le |g_n|{3e^{\|p\|_1}\/2}$. Using these estimates and $\b(z)=\b(e_n^-)+\b'(z_*)(z-e_n^-)$ for some $z_*\in g_n$ we obtain \er{T31-6}.

iii) Using \er{efs1} we obtain
$$
|(\D^2(z)-1)-(\cos^2z-1)|\le 2X|\D(z)-\cos z|\le 2X^2/|z|,\qq X=e^{|\Im |+\|p\|_1}.
$$
After this the standard arguments (due to Rouche's theorem) give
the proof of iii).
\BBox

 The function $a$ is not convenient, since $a$ is not analytic on $\cM$. The modified function $\x$ in the momentum variable $z$ is given by
\[
\lb{ix}
\x(z)=2i\sin k(z)a(z)=\Xi(z^2), \qqq z\in \cZ_0.
\]

The number $e_n=\sqrt{E_n}\in [e_n^-,e_n^+]$ satisfies $|\D(e_n)|=\max_{z\in g_n} |\D(z)|=\cosh h_n\ge 1$ for some $h_n\ge 0$.

\begin{lemma}
\lb{T32}
i) The following identity and asymptotics hold true:
\begin{multline}
\lb{T32-1}
\x(z)=2i\sin k(z)(1-A(z^2))-J(z^2),\qq z\in \cM,\qqq J(z^2)=\int_\R q(x)Y_1(x,z)dx
\\
A(z^2)=\int_\R q(x)Y_2(x,z)dx,\qqq
Y_1=\vp_1\vt \wt\vt-\vt_1'\vp \wt\vp+\b(\vp \wt\vt+\vt \wt\vp),
\qqq Y_2={1\/2}(\vp \wt\vt-\vt \wt\vp),\\
\x(z)=2(-1)^{n+1}(1+A(z^2))\sinh v(z)-J(z^2),\qqq z\in g_n^\pm\ne \es,
\end{multline}
where $v=\Im k$  and  $\pm v(z)>0$ for  $z\in g_n^\pm$,
\begin{multline}
\lb{T32-2}
Y_1=\vp(1,z,x)+Y_{11}(z,x),\ \
Y_{11}=\vp_1\vt \wt\vt_*-\vt_1'\vp \wt\vp_*+\b(\vp \wt\vt_*+\vt \wt\vp_*),
\ \
\wt\vt_*=\wt\vt-\vt, \wt\vp_*=\wt\vp-\vp,
\\
Y_{11}(z,x)=O(|g_n|/n^2)\ \ as \ \ z\in g_n, n\to \iy.
\end{multline}
Moreover, the functions $J, A$ are entire and $\x$ is analytic on $\cM$.

\no ii)  Let $|z|\ge 2$. Then the following estimates hold true
\begin{multline}
\lb{T32-3}
|J(z^2)|\le C_{p,q}\|q\|_t \rt(|\vp(1,z)|+{|\vt'(1,z)|\/|z|^2}+{|\b(z)|\/|z|}\rt)e^{2t|\Im z|}
\le {3C_*\/|z|}e^{(2t+1)|\Im z|},
\\
|A(z^2)|\le {\|q\|_t^2C_{p,q}\/|z|^2}e^{2t|\Im z|},\qq where \qq C_*=\|q\|_te^{\|p\|_1+\|q\|_t+2\|p\|_t}, \qq
C_{p,q}=e^{2\|p\|_t+\|q\|_t},
\end{multline}
\[
\lb{T32-4}
|J(e_n^2)|\le {R_1\/|z_n|}\sinh h_n,\qqq if \qq e_n^-\ge R_p=8e^{\|p\|_1}, \qqq R_1=11C_*,
\]
\end{lemma}
\no {\bf Proof.} i)
We rewrite the identity \er{T23-2} in the form
\[
\lb{318}
\x(z)=2i\sin k(z)-\int_\R q(x)Y(x,z)dx,\qq
Y=\vp_1\p_-(x,z)f_+(x,z),
\]
for $z\in \C_+$. Using \er{318}, \er{T22-1} we rewrite $Y$ in the form
$$
Y=\vp_1(\vt+m_-\vp)(\wt\vt+m_+\wt\vp)=\vp_1\rt(\vt \wt\vt
+m_+\vt \wt\vp+m_-\vp \wt\vt-{\vt_1'\/\vp_1}\vp \wt\vp\rt)
=Y_1-i2Y_2\sin k.
$$
Substituting the last identity into $\int_\R q(x)Y(x,z)dx$ and
using \er{35},  we obtain \er{T32-1}.

Substituting asymptotics from \er{efs}, Lemma \ref{T31} into $Y_{11}$  we obtain  $Y_{11}(z,x)=O(|g_n|/n^2)$ as $z\in g_n, n\to \iy.$

 ii) Using \er{efs}, \er{efs1}  and  \er{T32-1} and Lemma \ref{T31}, we obtain \er{T32-3}. Estimates \er{T32-3} and Lemma \ref{T31} give
$|J(e_n^2)|\le {21C_*\/4|z_n|} |g_n|$;
and the estimate $|g_n|\le 2|h_n|\le 2\sinh h_n$ from \cite{KK} yields
\er{T32-4}.
\BBox

This lemma gives that the function $\x$ is analytic on
$\cM$ and the function $\Xi$ is analytic on
$\L$. We define the bound states, resonances in terms of momentum variable
$z\in \cM$. Recall that there are bound states on the physical gaps and resonances on the non physical gaps.
Define the set $\cM_0 =\{z\in \cM: z=e_{\pm n}^+\in \cM_1 $ and $z=e_{\pm n}^+\in \cM_2, \g_n=\es, n\ge 1\}\ss \cM$. The set $\cM_0$ is the image of $\L_0$
(see before Definition of states, Section 1)
under the mapping $z=\sqrt \l$.

\no {\bf Definition.} {\it Each zero of $\x(z)=2i\sin k(z)a(z)$ in $\cM\sm \cM_0$ is a state  of $H$.

\no 1) A state  $z\in \cM_2$  is a resonance.

\no 2) A state   $z=e_n^\pm , n\in \Z$ is a virtual state.

\no 3) A resonance $z\in \cup g_n^{-}\ss\cM_2$ is an anti-bound state.}

Of course, $z^2$ is really the energy, but since
the momentum $z$ is the natural parameter, we will abuse the terminology.

The kernel of the resolvent $R=(H-z^2)^{-1},  z\in \C_+,$ has the form
$$
R(x,x',z)={f_- (x,z )f_+(x',z)\/-w(z)}={R_1(x,x',z)\/-\x(z)},\ \ \ x<x',\ \
R_1=\vp(1,z)f_- (x,z )f_+(x',z),
$$
and $R(x,x',z )=R(x',x,z ),\ x>x'$. Identity \er{T22-1} yields
$f_\pm=\wt\vt+m_\pm \wt\vp,\qq \wt\vt=\wt\vt(x),\qq  \wt\vp=\wt\vp(x)$. Let $ \wt \vt_*=\wt\vt(x'),\qq \wt \vp_*=\wt\vp(x')$. Then using \er{LD0} we obtain
$$
R_1(x,x',z)=\vp(1,\cdot)\wt \vt \wt \vt_*+(\b-i\sin k)\wt \vp \wt \vt_*+
(\b+i\sin k)\wt\vt \wt \vp_*-\vt'(1,\cdot)\wt \vp_*\wt \vp.
$$
Thus if
$\x(z)=\vp(1,z)w(z)=0$ at some $z\in\cM$, then  $(H-z^2)^{-1}$ has singularity at $z$. The poles of $R(x,x',z)$ define the bound states and resonances. The zeros of $\x$ define the bound states and resonances, since the function $R_1=\vp_1(z)f^-(x,z )f^+(x',z)$ is locally bounded.

If $q=0$, then $R_0=(H_0-z^2)^{-1}$ has the form
$$
R_0(x,x',z)={R_{10}(x,x',z)\/-\x_0(z)},\qq \x_0=\vp(1,z)w_0(z)=2i\sin k(z),
$$$$
R_{10}=\vp_1(z)\p_- (x,z )\p_+(x',z)=
\vp(1,\cdot)\vt\vt_*+(\b-i\sin k)\vp\vt_*+
(\b+i\sin k)\vt\vp_*-\vt'(1,\cdot)\vp\vp_*,
$$
where $\vp=\vp(x,z ),..$ and $\vp_*=\vp(x',z ),..$ Thus $R_0(x,x',z)$ has singularity
at some $z\in\cM$ iff $\sin k(z)=0$, i.e., $k(z)=\pi n $ and then
$z=e_n^\pm$.

Define the functions $F, S$ by
$$
F=\x(z)\x(-z), \qqq S(z)=\vp^2(1,z)s(z)s(-z), \qq z\in\cZ_0.
$$

\begin{lemma}
\lb{T33}
i)  The functions $F(z), S(z), z\in\cZ_0$
have analytic continuations into the whole complex plane $\C$ and satisfy
\[
\lb{T33-1}
F(z)=4(1-\D^2(z))(1+A(z^2))^2+J^2(z^2)=4(1-\D^2(z^2))+S(z), \qqq z\in \C.
\]
Moreover, $F(z)>0$ and $S(z)\ge 0$ on each interval  $(e_{n-1}^+,e_n^-), n\ge 1$ and $F$ has even number of zeros on each interval $[e_n^-,e_n^+], n\ge 1$. The function $F$ has only simple zeros at $e_n^\pm, g_n\ne \es$.

\no ii) If $g_n=(e_n^-,e_n^+)= \es$ for some $n\ne 0$, then each
$f_\pm(x,\cdot ),  x\in \R$ is analytic in some disk $\cD(e_n^+,\ve),\ve >0$. Moreover, $\m_n=e_n^\pm$ is
a double zero of $F$ and $e_n^+$ is not a state of $H$.

\no iii) Let $z\in g_n^+$ be a bound state for some $n\ge 1$, i.e.,
$\x(z+i0)=0$. Then $z-i0\in g_n^-$ is not an anti-bound state
and $\x(z-i0)\ne 0$.

\no iv) Let $z\in i\R_+$ be a bound state, i.e.,
$\x(z)=0$. Then $-z\in i\R_-$ is not an anti-bound state and $\x(-z)\ne 0$.

\no v) $z\in \C_-\sm i\R$ is a zero of $F$ iff $z\in \C_-\sm i\R$ is a zero of $\x$ (with the same multiplicity).

\no vi) $z\in i\R_-$ is a zero of $F$ iff $z\in i\R_-$ or
$-z\in i\R_-$ is a zero of $\x$.

\no vii) Let $g_n\ne \es, n\ge 1$. The point $z\in g_n$ is a zero of $F$ iff $z+i0\in g_n^+$ or $z-i0\in g_n^-$
 is a zero of $\x$ (with the same multiplicity).

\no viii) Let $e_n^-> \max\{8e^{\|p\|_1}, 11C_*\}, C_*=\|q\|_te^{\|p\|_1+\|q\|_t+2\|p\|_t}$ and $|g_n|>0$. Then
\[
\lb{T33-2}
F(e_n)\le - \rt(1-{R_1^2\/|e_n|^2}  \rt)\sinh^2 h_n<0, \qqq \where
\qqq e_n=\sqrt {E_n}>0,
\]
and $E_n,h_n$ are defined by \er{dzn}. Moreover, $F$  has at list two zeros in the segment $[e_n^-,e_n^+]$.

\end{lemma}

\no{\bf Proof.} i) Using \er{T32-1} we deduce that
$F=\x(z)\x(-z)$ is entire and satisfies $F(z)=(1-\D^2(z))(2-A(z^2))^2+J^2(z^2)$.
Using \er{iab} we obtain that $S$ is entire  and satisfies \er{T33-1}.

Recall that $F\ge 0$ and $S\ge0$ inside each $(e_{n-1}^+,e_n^-)$.
Moreover, $F>0$, since $a\ne 0$ inside each $(e_{n-1}^+,e_n^-)$.
Due to $F(e_n^\pm)\ge 0$, we get that $F$ has even number of zeros on each interval $[e_n^-,e_n^+], n\ge 1$.
We have $F=F_0+S, F_0=4(1-\D^2)$.
Consider the case $z=e_n^+$, the proof for $z=e_n^-$ is similar.
Thus if $F(z)=0$, then we get $S(z)=0$, since
$\D^2(e_n^+)=1$. Moreover, $F_0'(e_n^+)=-2\D(e_n^+)\D'(e_n^+)>0$
and $S'(e_n^+)\ge 0$, which gives that $z=e_n^+$ is a simple zero of $F$.

ii) The function $m_\pm$ is analytic in $\cZ_0$, then each $f^\pm(x,\cdot), x\in \R$ is analytic in $\cZ_0$.
Moreover, \er{LD0} yields $\b(\m_n)=\vp(1,\m_n)=\vt'(1,\m_n)=0$ and then $J(\m_n^2)=0$. Thus by \er{T33-1}, the function $F$ has a double zero at $\m_n$ at list.

 iii)  If $z\in g_n, n\ge 1$ is a bound state. Then
\er{T32-1} yields
$$
0=\x(z+i0)=2(-1)^{n+1}\sinh v(z+i0)(1+A(z^2))-J(z^2),
$$
$$
2\sinh v(z+i0)={(-1)^{n+1}J(z^2)|\/(1+A(z^2))}>0.
$$
If we assume that $z-i0$ is a anti-bound state, then $2\sinh v(z-i0)={(-1)^{n+1}J(z^2)\/(1+A(z^2))}<0$, which gives contradiction. The proof of iv) is similar.

v) The function $\x$ has not zeros in $\C_+\sm iR$, see \cite{F1}.
This yields v).

vi) The property $F(z)=F(-z), z\in \C$ gives vi).

vii)  The statement vii) follows from iii).

viii) Estimate \er{T32-3} gives $|A(e_n^2)|\le {1\/2}$ and \er{T32-4} imply \er{T33-2}. The function $F(e_n^\pm)\ge 0$ and due to \er{T33-2}, we deduce that $F$  has at list two zeros in the segment $[e_n^-,e_n^+]$.
\BBox

{\bf Proof of Theorem \ref{T1}}. Identity \er{T1-1} have been proved in Lemma
\ref{T33}.

i) The function $\x$ is analytic in
$\cZ_0$ and is real on $i\R$. Then the set of zeros of $\x$ is
symmetric with respect to the imaginary line and satisfies
\er{T1-1}.

The statement ii) has been proved in Lemma \ref{T33}.

iii) Using Lemma \ref{T32}, \er{efs} and \er{efs1}  we obtain
$$
|A(z^2)|\le {\|q\|_t^2\/|z|^2}X_2,\qqq |J(z^2)|\le {3\|q\|_t\/|z|}
X_2X\qqq |\D(z)-\cos z|\le {X\/|z|},\qq |z|\ge 2,
$$
where $ X_2=e^{|\Im z|2t+\|q\|_t+2\|p\|_t},\qq X=e^{|\Im z|+\|p\|_1} $. Substituting these estimates into the identity
$$
F-4\sin^2 z=4(\cos^2z-\D^2)+J^2+4(1-\D^2)A(A+2)
$$
we obtain
$$
|F(z)-4\sin^2 z|\le 9X^2C_0,\qqq
C_0={1\/|z|}+{\|q\|_t^2X_2^2\/|z|^2}+{\|q\|_t^2X_2\/|z|^2}\rt(2+{\|q\|_t^2X_2\/|z|^2}\rt),\qq |z|\ge 2.
$$
Using $e^{|{\Im}z|}\le 4|\sin z|$
for all $|z-\pi n|\ge {\pi \/4}, n\in \Z$, (see p. 27 [PT]), we obtain
$$
9X^2=9e^{2|\Im z|+2\|p\|_1}\le  |4\sin^2 z|r_0^2
 \qq all \ |z-\pi n|\ge{\pi\/4},\qq  n\in \Z, \qqq r_0=6e^{\|p\|_1},
$$
which yields
$$
|F(z)-4\sin^2 z|\le 4|\sin^2 z|{r_0^2C_0 \/|z|}< 4|\sin^2 z|,\qq all
\qq |z|\ge 2, z\notin \cup \cD(\pi n,{\pi\/4}),
$$
since
$$
{r_0^2C_0\/|z|}\le
 {r_0^2\/|z|}+{r_0^2\|q\|_t^2X_2^2\/|z|^2}+{r_0^2\|q\|_t^2X_2\/|z|^2}\rt(2+{\|q\|_t^2X_2\/|z|^2}\rt)<{1\/5}
+{1\/4}+{1\/2}+{1\/(24)^2}<{19\/20},
$$
where  we have used:
$$
if \ z\in \cD_F=\{z\in \C: |z|>\max \{180e^{2\|p\|_1},C_Fe^{2t|\Im z|}\}\}\qq \Rightarrow \qq
{r_0^2\/|z|}<{1\/5},\qq {r_0\|q\|_tX_2\/|z|}<{1\/2},
$$
and recall that
$$
\mD_F=\{\l\in \L_2: |\l|^{1\/2}>\max \{180e^{2\|p\|_1},C_Fe^{2t|\Im \sqrt \l|}\}\},\ \ C_F=12\|q\|_te^{\|p\|_1+\|q\|_t+2\|p\|_t}.
$$
Thus  by Rouche's theorem, $F$ has as many roots, counted with
multiplicities, as $\sin^2 z$ in each disk $\cD_{\pi\/4}(\pi n)\ss\cD_F$.
Since $\sin z$ has only the roots $\pi n, n\ge 1$, then $F$
has two zeros in each disk $\cD_{\pi\/4}(\pi n)\ss\cD_F$ and $F$ has not zeros in $\cD_F\sm \cup \cD_{\pi\/4}(\pi n)$.

$F(e_n^\pm)\ge 0$ for all $n$ and Lemma \ref{T33} yields $F(e_n)<0$ for all $n$,
where $e_n^-> \max\{8e^{\|p\|_1}, 11C_*\}$, see \er{T33-2}. Then
there exists exactly two   simple real zeros  $\vk_n^\pm=\sqrt{\l_n^\pm}>0$ of $F$
such that  $e_n^-\le \vk_n^-<e_n<\vk_n^+\le e_n^+$. Moreover, $F$
has not no states in the "logarithmic"  domain $\cD_F\cap \C_-$.

 The statements iv) and v) have been proved in Lemma \ref{T33}.
\BBox

\section {Proof of Theorem 1.2-1.4}
\setcounter{equation}{0}

\no {\bf Proof of Theorem \ref{T2}}. i) Theorem \ref{T1} gives
 $\vk_n^\pm=e_n^\pm\mp \d_n^\pm=\sqrt{\l_n^\pm}$. Let $\vk=\vk_n^\pm,  \d=\d_n^\pm$.
Then the equation $0=\x(\vk)=(-1)^{n+1}2(1+A(\vk^2))\sinh v(\vk)-J(\vk^2),\  \ \vk\in \ol g_n^\pm\ne \es$ and \er{pav} imply
$$
\sinh |v(\vk)|=O(J(\vk^2))=\ve O( |\vp(1,\vk)|+{|\vt'(1,\vk)|\/|\vk|^2}+{|\b(\vk)|\/\vk}) =\ve O(|g_n|), \qq \ve={1\/2\pi n}
$$
as $n\to \iy$.
 Moreover, using the estimate $|(z-e_n^-)(z-e_n^+)|^{1\/2}\le |v(z)|$ for each $z\in g_n$ (see \cite{KK}) we obtain
$|\d(|g_n|-\d)|^{1\/2}\le |v(\vk)|=\ve O(|g_n|)$,
which yields $\d=\ve^2 O(|g_n|)$. Thus the points $\vk_n^\pm$ are close to $e_n^\pm$ and satisfy:
\[
\lb{asd}
\d_n^-=\vk_n^--e_n^-=\ve^2 O(|g_n|),\qqq and \qqq
\d_n^+=e_n^+-\vk_n^+=\ve^2 O(|g_n|).
\]
Recall that $\vk_n^-$ and $\vk_n^+$ are simple.
Consider the first case $\d=\d_n^-=|g_n|O(\ve^2)$.

Using \er{T32-2} we obtain
$$
J=J_{10}+J_{10},\qqq J_{10}(z)=\int_\R\vp(1,z,x)q(x)dx.
$$
Let $\m_n^2(\t), \t\in \R$ be the Dirichlet eigenvalue for the problem $-y''+q(x+\t)y=z^2 y, y(0)=y(1)=0$.
In this case $\cos \f_n$ in \er{ip1} is a function from $\t\in \R$.
Below we need facts from \cite{K5}:
\[
\lb{mz}
{E_n^-+E_n^+\/2}-\m_n^2(\t)={|\g_n|\/2}\cos \f_n(\t),\qq
\f_n(\t)=\f_n(0)+2\pi n\t+O(\ve),
\]
as $n\to \iy$ uniformly with respect to $\t\in [0,1]$.
Asymptotics \er{T31-1} yield
\[
\vp(1,\vk,x)={(-1)^n\/\pi n}(1+O(\ve))(\vk-\m_n(x))={(-1)^n\/\pi n}(1+O(\ve))(e_n^--\m_n(x)+\ve^2 O(|g_n|)).
\]
Thus we rewrite $e_n^--\m_n(x)$ in the form
$$
e_n^--\m_n(x)={{E_n^-+E_n^+-|\g_n|\/2}-\m_n^2(x)\/e_n^-+\m_n(x)}=
{|\g_n|\/2} {\cos \f_n(x)-1\/e_n^-+\m_n(x)}.
$$
This gives
$$
J_{10}(\vk)=\int_\R\vp(1,\vk,x)q(x)dx={(-1)^n\/\pi n}(1+O(\ve))
\int_\R \rt[{|\g_n|\/2} {\cos \f_n(x)-1\/e_n^-+\m_n(x)}+\ve^2 O(|g_n|)\rt]     q(x)dx.
$$
$$
={(-1)^n|g_n|\/2\pi n}
\int_\R \rt[\cos y_n(x)-1+ O(\ve)\rt]q(x)dx
={(-1)^n|g_n|\/2\pi n}(-q_0+c_n\wh q_{cn}-s_n\wh q_{sn}+ O(\ve)).
$$
where $c_n=\cos y_n(0), s_n=\sin y_n(0)$, and thus
\[
(-1)^{n+1}J(\vk^2)={|g_n|\/2\pi n}I_n^-, \qqq I_n^-=q_0-c_n\wh q_{cn}+s_n\wh q_{sn}+ O(\ve).
\]
Using \er{T32-3} we obtain
\[
\lb{as2}
\sinh v(\vk)={(-1)^{n+1}J(\vk^2)\/2+2A(\vk^2)}={\ve |g_n|I_n^-\/2+O(\ve^2)},\qq
\sign v(\vk)=\sign (-1)^{n+1}J(\vk^2)=\sign I_n^-.
\]
Note  that if $v(\vk)>0$ (or $v(\vk)<0$), then $\vk\in g_n^+$ (or $\vk\in g_n^-$) and $\vk$ is a bound state (or a resonance). Moreover,  if $v(\vk)=0$, then  $\vk=e_n^-$ or $\vk=e_n^+$ is a virtual state. Then  \er{asd} gives
$\sinh v(\vk)=v(\vk)(1+O(|g_n|^2\ve^2))$  and using asymptotics \er{pav}, we obtain
$$
 v(\vk)=\sqrt{\d(|g_n|-\d)}(1+O(\ve^2))=\sqrt{\d|g_n|}(1+O(\ve^2)).
$$
Thus \er{as2} yields $\d_n^-={|g_n|\ve^2\/4}(I_n^-)^2$ and \er{ape} gives $|\g_n|=(2\pi n)|g_n|(1+O(\ve^2))$,
which yields $\d_n^-={2|\g_n|\/(4\pi n)^3}(I_n^-)^2(1+O(\ve^2))$.

If $q_0>0$, then $I_n^->0$ and  above arguments yield that $\vk_n^-$ is a bound state and
$\vk_n^+$ is an anti bound state. Conversely, if $q_0<0$, then $I_n^-<0$ and we deduce that $\vk_n^-$ is an anti  bound state and $\vk_n^+$ is a bound state.

Similar  arguments imply the proof for the case $\d_n^+=e_n^+-\vk_n^+=|g_n|O(\ve^2)$.

iii) We have $-c_n\wh q_{cn}+s_n\wh q_{sn}=-|\wh q_n|\cos (\f_n+\t_n)$. Then (i) yields
the statement (ii).

iv) An entire function $f(z)$ is said to be of $exponential$ $ type$  if
there is a constant $\a$ such that  $|f(z)|\leq $ const. $e^{\a |z|}$
everywhere. The function $f$ is said to belong to the Cartwright class $\mE_\r,$
if $f$ is entire, of exponential type, and the following conditions are fulfilled:
$$
\int _{\R}{\log ^+|f(x)|dx\/ 1+x^2}<\iy  ,\ \
\r_\pm(f)=\r,\ \ \ {\rm where}\ \ \
\r_{\pm}(f)\ev \lim \sup_{y\to \iy} {\log |f(\pm iy)|\/y}.
$$

 Denote by $\cN^+(r,f)$ the number of zeros of $f$ with real
part $\geq 0$ having modulus $\leq r$, and by $\cN^-(r,f)$ the
number of its zeros with real part $< 0$ having modulus $\leq r$,
each zero being counted according to its multiplicity. We recall the
well known result (see [Koo]).

\no    {\bf Theorem (Levinson).}
{\it  Let the entire function $f\in \mE_\r,\r>0$.
Then $  \cN^{\pm }(r,f)={r\/ \pi }(\r+o(1))$ as $r\to \iy$,
and for each $\d >0$ the number of zeros of $f$ with modulus $\leq r$
lying outside both of the two sectors $|\arg z | , |\arg z -\pi |<\d$
is $o(r)$ for $r\to \iy$.}

Let $\cN (r,f)$ be the total number of  zeros of $f$  with modulus $\le r$.
Denote by $\cN_+(r,f)$ (or $\cN_-(r,f)$) the number of zeros of $f$ with imaginary part $>0$ (or $<0$)
 having modulus $\le r$, each zero being counted according to its multiplicity.

Let $\pm \z_n>0, n\in \N$  be all real zeros $\ne 0$  of $F$
and let the zero $\z_0=0$ has the multiplicity $n_0\le 2$.
Let $F_1=z^{n_0}\lim_{r\to \iy}\prod_{|\z_n|\le r}(1-{z\/\z_n})$.
The Levinson Theorem and Lemma \ref{T23} imply
$$
\cN(r,F)=\cN(r,F_1)+\cN(r,F/F_1)=2r{1+2t+o(1)\/\pi},\qq
\cN(r,F_1)=2r{1+o(1)\/\pi} \ \ \as \ r\to\iy.
$$
Then Lemma \ref{T33} gives the identities  $\cN_-(r,F)=\cN_+(r,F)=\cN_-(r,\x)+N_*$ for some integer $N_*\ge 0$. Thus we obtain
\[
\cN(r,F)=\cN(r,F_1)+2\cN_-(r,\x)+2N_*=2r{1+2t+o(1)\/\pi},
\]
which yields $\cN_-(r,\x)={2rt+o(r)\/\pi}$ as $r\to \iy$ and \er{T2-2}.

Due to \er{T2-1}, the high energy states of $H$ and $H_0$ are very close.
 This gives \er{T2-3}.
 \BBox

\no {\bf Proof of Theorem \ref{T3}}. i)
 Let the operator $H_0$ have infinitely many gaps $\g_n\ne \es$
for some $p\in L^2(0,1)$ and let $\s=(\s_n)_{1}^\iy$ be any sequence,
where $\s_n\in \{0,2\}$.

We take a potential $q$ and let $\wh q_n=\wh q_{cn}+i\wh q_{sn}=|\wh q_n|e^{i\t_n}, n\ge 1$ and $q_0=0$.
We also assume that $|\wh q_n|>n^{-\a}$ for  $n$ large enough and for some  $\a\in (0,1)$.
For each $\f_n$ we take $\t_n$ such that $|\cos (\f_n+\t_n)|>\ve>0$. Then due to Theorem \ref{T2} iii), the operator $H$ has $\s_n=1-\sign \cos (\f_n+\t_n)$ bound states in the physical gap $\g_n^{(1)}\ne \es$
and $2-\s_n$ resonances inside the  nonphysical gap $\g_n^{(2)}\ne \es$  for $n$ large enough.
Thus changing $\t_n$ we obtain $\cos (\f_n+\t_n)<0$, which yields  $\s_n=2$ or we obtain $\sign \cos (\f_n+\t_n)>0$, which yields $\s_n=0$.

ii) Let  $q\in \cQ_t,t>0$ satisfy   $q_0=0$ and let $|\wh q_{cn}|>n^{-\a}$ for all $n$ large enough and some $\a\in (0,1)$. Let $\s=(\s_n)_{1}^\iy$ be any sequence, where $\s_n\in \{0,2\}$.
Let $\d=(\d_n)_1^\iy\in \ell^2$ be a sequence of nonnegative numbers $\d_n\ge 0, n\ge 1$ and
infinitely many $\d_n>0$.

Recall the result from \cite{K5}:

{\it The mapping $\P: \cH\to \ell^2\os \ell^2$ given by $\P=((\P_{cn})_1^\iy,(\P_{sn})_1^\iy)$ is a real analytic isomorphism between real Hilbert spaces $\cH=\{p\in L^2(0,1): \int_0^1p(x)dx=0\}$ and $\ell^2\os \ell^2$.}

Then for the sequence $\d=(\d_n)_1^\iy\in \ell^2$ there exists a potential $p\in L^2(0,1)$
such that each gap length  $|\g_n|=\d_n$. Assume that $E_0^+=0$. Moreover, using Theorem \ref{T2} we deduce that
$H$ has exactly two  simple states $\l_n^-,\l_n^+\in \g_n^{(0)}$ for $n$ large enough and $\l_n^-,\l_n^+$
have asymptotics
$$
\sqrt{\l_n^\pm}=\sqrt{E_n^\pm}\mp {2|\g_n|\/(4\pi n)^3}(R_n+ O(1/n))^2,\qq
\qqq (-1)^{n+1}J(\l_n^\pm)={|\g_n|\/(2\pi n)^2}(R_n+ O(1/n)),
$$
where $R_n=-c_n\wh q_{cn}+s_n\wh q_{sn}=-|\wh q_n|\cos (\f_n+\t_n)$ and
$\wh q_n=\wh q_{cn}+i\wh q_{sn}=|\wh q_n|e^{i\t_n}$ and $\f_n$ is defined in \er{ip1}.
The parameters $\t_n$ are fixed, but due to above results from \cite{K5} the angles $\f_n$ can be any numbers.
If we take $\f_n$ such that $|R_n|>\ve>0$ for all $\g_n\ne \es$
and some $\ve>0$. Moreover, we take $\f_n$ such that $\s_n=1-\sign \cos (\f_n+\t_n)$.
Then Theorem \ref{T2} we obtain that the operator $H$ has
$\s_n=1$ bound states in the physical gap $\g_n^{(1)}\ne \es$
and $2-\s_n$ resonances inside the  nonphysical gap $\g_n^{(2)}\ne \es$  for $n$ large enough.
\BBox

\no {\bf Acknowledgments.}
\small
The various parts of this paper were written at  ESI, Vienna,
and Mathematical Institute of the Tsukuba Univ., Japan and Ecole Polytechnique, France.
The author is grateful to the Institutes for the hospitality.

%\bigskip

\end{document}